\documentclass[fleqn]{mat01}
\usepackage{times,mathtimy,amssymb,latexsym,amsbsy}
\begin{document}

\setcounter{page}{79}
\firstpage{79}

\newtheorem{theore}[subsubsection]{Theorem}
\newtheorem{theor}[subsubsection]{\bf Theorem}
\newtheorem{lem}[subsubsection]{Lemma}
\newtheorem{rem}[subsubsection]{Remark}
\newtheorem{pot}[subsubsection]{Proof of Theorem}
\newtheorem{coro}[subsubsection]{\rm COROLLARY}

\newtheorem{theo}[ssection]{\bf Theorem}

\newtheorem{lemm}{Lemma}
\renewcommand\thelemm{\it \Alph{lemm}}

\title{Inequalities for dual quermassintegrals of mixed intersection bodies}

\markboth{Zhao Chang-jian and Leng Gang-song}{Inequalities for dual quermassintegrals}

\author{ZHAO CHANG-JIAN$^{1}$ and LENG GANG-SONG$^{2}$}

\address{$^{1}$Department of Information and Mathematics Sciences,
College of Science, China~Institute of Metrology, Hangzhou~310018,
People's Republic of China\\
\noindent $^{2}$Department of Mathematics, Shanghai University, Shanghai 200~436,
People's~Republic of China\\
\noindent E-mail: chjzhao@163.com; lenggangsong@163.com}

\volume{115}

\mon{February}

\parts{1}

\pubyear{2005}

\Date{MS received 8 May 2004}

\begin{abstract}
In this paper, we first introduce a new concept
of {\it dual quermassintegral sum function} of two star bodies and
establish Minkowski's type inequality for dual quermassintegral
sum of mixed intersection bodies, which is a general form of the
Minkowski inequality for mixed intersection bodies. Then, we give
the Aleksandrov--Fenchel inequality and the Brunn--Minkowski
inequality for mixed intersection bodies and some related results.
Our results present, for intersection bodies, all dual
inequalities for Lutwak's mixed prosection bodies inequalities.
\end{abstract}


\maketitle

\setcounter{section}{-1}
\section{Introduction}

One might say the history of intersection bodies began with the
paper of Busemann \cite{4}. Intersection bodies were first explicitly
defined and named by Lutwak \cite{11}. It was here that the duality
between intersection bodies and projection bodies was first made
clear. Despite considerable ingenuity of earlier attacks on
the Busemann--Petty problem, it seems fair to say that the work
of Lutwak \cite {11} represents the beginning of its eventual solution. In
\cite{11}, Lutwak also showed that if a convex body is sufficiently
smooth and not an intersection body, then there exists a centred
star body such that the conditions of Busemann--Petty problem
holds, but the result inequality is reversed. Following Lutwak,
the intersection body of order $i$ of a star body is introduced by
Zhang \cite{21}. It follows from this definition that every intersection
body of order $i$ of a star body is an intersection body of a star
body, and vice versa. As Zhang observes, the new definition of
intersection body allows a more appealing formulation, namely: the
Busemann--Petty problem has a positive answer in $n$-dimensional
Euclidean space if and only if each centered convex body is an
intersection body. The intersection body plays an essential role
in Busemann's theory \cite{5} of area in Minkowski spaces. The
intersection body is also an important matter of the
Brunn--Minkowski theory.

In recent years, some authors including Ball \cite{1,2}, Bourgain \cite{3},
Gardner \cite{6,7,8}, Schneider \cite{19} and Lutwak \cite{12,13,14,15,16,17,18} have given
considerable attention to the Brunn--Minkowski theory and their
various generalizations. The purpose of this paper is to
establish the Minkowski inequality for the dual quermassintegral
sum, which is a generalization of the Minkowski inequality for
mixed intersection bodies. Then, the Brunn--Minkowski inequality
and the Aleksandrov--Fenchel inequality for mixed intersection
bodies are proved and some related results are also given. In this
work we shall derive, for intersection bodies, all the analogous
inequalities for Lutwak's mixed projection body inequalities \cite{15}.
Thus, this work may be seen as presenting additional evidence of
the natural duality between intersection and projection bodies.

\section{Notation and preliminaries}

The setting for this paper is an $n$-dimensional Euclidean space
${\Bbb R}^{n} \ (n>2)$. Let ${\Bbb C}_{n}$ denote the set of
non-empty convex figures (compact, convex subsets) and let ${\cal
K}^{n}$ denote the subset of  ${\Bbb C}_{n}$ consisting of all
convex bodies (compact, convex subsets with non-empty interiors)
in ${\Bbb R}^n$. We reserve $u$ for unit vectors, and
$B$ for the unit ball centered at the
origin. The surface of $B$ is $S^{n-1}$. For $u\in S^{n-1}$, let
$E_u$ denote the hyperplane, through the origin, that is
orthogonal to $u$. Let $K^{u}$ to denote the image of $K$
under an orthogonal projection onto the hyperplane $E_u$. We use
$V(K)$ for the $n$-dimensional volume of convex body $K$. The
support function of $K\in {\cal K}^n$, $h(K,\cdot)$, defined on
${\Bbb R}^{n}$ by $h(K,\cdot)=\max\{x\cdot y: y\in K\}.$ Let
$\delta$ denote the Hausdorff metric on ${\cal K}^n$; i.e., for
$K, L\in {\cal K}^{n}$, $\delta(K,L)=|h_K-h_L|_{\infty}$, where
$|\cdot|_{\infty}$ denotes the sup-norm on the space of continuous
functions, $C(S^{n-1}).$

Associated with a compact subset $K$ of ${\Bbb R}^n$, which is
star-shaped with respect to the origin, its radial function
$\rho(K,\cdot): S^{n-1}\rightarrow {\Bbb R}$, defined for $u\in
S^{n-1}$, by $\rho(K,u)=\max\{\lambda\geq 0: \lambda u\in K\}$. If
$\rho(K,\cdot)$ is positive and continuous, $K$ will be called a
star body. Let $\varphi^{n}$ denote the set of star bodies in
${\Bbb R}^{n}$.\vspace{.5pc}

\subsection{\it Dual mixed volumes}

If $K_{1},\ldots,K_{r}\in \varphi^n$ and
$\lambda_{1},\ldots,\lambda_{r}\in {\Bbb R}$, then the radial
Minkowski  linear combination,
$\lambda_{1}K_{1}\tilde{+}\cdots\tilde{+}\lambda_{r} K_{r}$, is
defined by $\lambda_{1}K_{1}\tilde{+}\cdots\tilde{+}\lambda_{r}
K_{r}=\{\lambda_{1}x_{1}\tilde{+}\cdots\tilde{+}\lambda_{r} x_{r}$:
$x_{i}\in K_{i}\}$.

The following property will be used later. If $K, L\in \varphi^n$
and $\lambda, \mu\geq 0$,
\renewcommand\theequation{\thesubsection.\arabic{equation}}
\begin{equation}
\rho(\lambda K\tilde{+}\mu L,\cdot)=\lambda
\rho(K,\cdot)+\mu\rho(L,\cdot).
\end{equation}
For $K_{1},\ldots,K_{r}\in \varphi^n$ and
$\lambda_{1},\ldots,\lambda_{r}\geq 0$, the volume of the radial
Minkowski linear combination
$\lambda_{1}K_{1}\tilde{+}\cdots\tilde{+}\lambda_{r}K_{r}$ is a
homogeneous $n$th-degree polynomial in the $\lambda_{i}$ \cite{19},
\begin{equation}
V(\lambda_{1}K_{1}\tilde{+}\cdots\tilde{+}\lambda_{r}K_{r})=\sum \tilde{V}_{i_{1},\ldots,i_{n}}\lambda_{i_{1}}\cdots\lambda_{i_{n}},
\end{equation}
where the sum is taken over all $n$-tuples $(i_{1},\ldots,i_{n})$
whose entries are positive integers not exceeding $r$. If we
require the coefficients of the polynomial in (1.1.2) to be
symmetric in their arguments, then they are uniquely determined.
The coefficient $\tilde{V}_{i_{1},\ldots,i_{n}}$ is non-negative
and depends only on the bodies $K_{i_{1}},\ldots,K_{i_{n}}$. It is
written as $\tilde{V}(K_{i_{1}},\ldots,K_{i_{n}})$ and is called
the {\it dual mixed volume} of $K_{i_{1}},\ldots,K_{i_{n}}$. If
$K_{1}=\cdots=K_{n-i}=K$, $K_{n-i+1}=\cdots=K_{n}=L$, the dual
mixed volumes is written as $\tilde{V}_{i}(K,L)$ and the dual mixed
volumes $\tilde{V}_{i}(K,B)$ is written as $\tilde{W}_{i}(K)$.

For $K,L\in \varphi^n$ and $i\in R$, the $i$th dual mixed volume
of $K$ and $L$, $\tilde{V}_{i}(K,L)$, is defined by \cite{14}
\begin{equation*}
\tilde{V}_{i}(K,L)=\frac{1}{n}\int_{S^{n-1}}\rho(K,u)^{n-i}\rho(L,u)^{i}\hbox{d}S(u).
\end{equation*}
From the above identity, if $K\in \varphi^{n}$, $i\in {\Bbb R}$, then
\begin{equation}
\tilde{W}_{i}(K)=\frac{1}{n}\int_{S^{n-1}}\rho(K,u)^{n-i}\hbox{d}S(u).
\end{equation}

If $K_i\in \varphi^{n}\, (i=1,2,\ldots,n-1)$, then the dual mixed
volume of $K_{i}\cap E_{u}\, (i=1,2,\ldots,n-1)$ will be denoted by
$\tilde{v}(K_{1}\cap E_{u},\ldots,K_{n-1}\cap E_u)$. If
$K_{1}=\cdots=K_{n-1-i}=K$ and $K_{n-i}=\cdots=K_{n-1}=L,$ then
$\tilde{v}(K_{1}\cap E_{u},\ldots,K_{n-1}\cap E_{u})$ is written as
$\tilde{v}_{i}(K\cap E_{u},L\cap E_{u})$. If $L=B$, then
$\tilde{v}_{i}(K\cap E_ {u},B\cap E_{u})$ is written as
$\tilde{w}_{i}(K\cap E_{u})$.\vspace{.5pc}

\subsection{\it Intersection bodies}

For $K\in \varphi^n$, there is a unique star body $I\!K$ whose
radial function satisfies for $u\in S^{n-1}$,
\setcounter{equation}{0}
\begin{equation}
\rho(I\!K,u)=v(K\cap E_{u}).
\end{equation}
It is called the {\it intersection bodies} of $K$. From a result
of Busemann, it follows that $I\!K$ is a convex if $K$ is convex and
centrally symmetric with respect to the origin. Clearly any
intersection body is centred.

The volume of the intersection bodies is given by
$V(I\!K)=\frac{1}{n}\int_{S^{n-1}}v(K\cap E_{u})^{n}\hbox{d}S(u)$.

The mixed intersection bodies of $K_{1},\ldots,K_{n-1}\in
\varphi^{n}$, $I(K_{1},\ldots,K_{n-1})$, whose radial function is
defined by
\begin{equation}
\rho(I(K_{1},\ldots,K_{n-1}),u)=\tilde{v}(K_{1}\cap E_{u},\ldots,K_{n-1}\cap E_{u}),
\end{equation}
where $\tilde{v}$ is $(n-1)$-dimensional dual mixed volume.

If $K\in \varphi^{n}$ with $\rho(K,u)\in C(S^{n-1})$, and $i\in
{\Bbb R}$ is positive, the {\it intersection body of order i} of
$K$ is the centered star body $I_{i}K$ such that \cite{3}
$\rho(I_{i}K)=\frac{1}{n-1}\int_{S^{n-1}}\rho(K,u)^{n-i-1}\hbox{d}S(u)$,
for $u\in S^{n-1}$, where
$I_{i}K=I(\underbrace{K,\ldots,K}_{n-i-1},\underbrace{B,\ldots,B}_{i})$.

If $K_{1}=\cdots=K_{n-i-1}=K, K_{n-i}=\cdots=K_{n-1}=L$, then
$I(K_{1},\ldots,K_{n-1})$ is written as $I_{i}(K,L)$. If $L=B$,
then $I_{i}(K,L)$ is written as $I_{i}K$ and is called the $i$th
intersection body of $K$. For $I_{0}K$ simply write $I\!K$. The term
is introduced by Zhang \cite{21}.

The following properties will be used later: If $K, L$, $M,
K_{1},\ldots,K_{n-1}\in \varphi^{n}$ and $\lambda, \mu,
\lambda_{1},\ldots,\lambda_{n-1}>0$, then
\begin{align}
&I(\lambda K\tilde{+}\mu L, M) =\lambda I(K,M)\tilde{+}\mu I(L,M),\quad \hbox{where}\ M=(K_{1},\ldots,K_{n-2}).\\
&I(\lambda_{1}K_{1},\ldots,\lambda_{n-1}K_{n-1}) =\lambda_{1}\cdots\lambda_{n-1}I(K_{1},\ldots,K_{n-1}).
\end{align}

\section{Main results}

The following results will be required to prove our main Theorems.

\begin{lemm}
If $K,L\in \varphi^{n}${\rm ,} $0\leq i<n$ and $0<j<n-1${\rm ,} then
\begin{align*}
\tilde{W}_{i}(I\!K) &=\frac{1}{n}\int_{S^{n-1}}v(K\cap E_{u})^{n-i}{\rm d}S(u),\\[.5pc]
\tilde{W}_{i}(I_{j}K) &=\frac{1}{n}\int_{S^{n-1}}\tilde{w}_{j}(K\cap E_{u})^{n-i}{\rm d}S(u),
\end{align*}
\begin{equation*}
\tilde{W}_{i}(I_{j}(K,L)) =\frac{1}{n}\int_{S^{n-1}}\tilde{v}_{j}(K\cap E_{u},L\cap E_{u})^{n-i}{\rm d}S(u).
\end{equation*}
\end{lemm}

To prove this, we use (1.1.3) in conjunction with the fact (1.2.2).

\begin{lemm}
{\rm \cite{14}}. If $K_{1},\ldots,K_{n}\in \varphi^{n}${\rm ,} then
\begin{equation*}
\tilde{V}(K_{1},\ldots,K_{n})^{r}\leq \prod_{j=1}^{r}\tilde{V}(\underbrace{K_{j},\ldots,K_{j}}_{r},K_{r+1},\ldots,K_{n})
\end{equation*}
with equality if and only if $K_{1},\ldots,K_{n}$ are
all dilations of each other.
\end{lemm}

We shall need the following trivial elementary inequality.

\begin{lemm}
If $a, b\geq 0$ and $c, d>0${\rm ,} then for $0<p<1${\rm ,}
\begin{equation*}
(a+b)^{p}(c+d)^{1-p}\geq a^{p}c^{p-1}+b^{p}d^{p-1},
\end{equation*}
with equality if and only if $ad=bc$.
\end{lemm}

\begin{proof}
Consideration the following function
\begin{equation*}
f(x)=(x+b)^{p}(c+d)^{1-p}-x^{p}c^{1-p},\qquad x\geq 0.
\end{equation*}

Let $f'(x)=p(c+d)^{1-p}(x+b)^{p-1}-p c^{1-p}x^{p-1}=0$, we get
$x={bc}/{d}$. If $x\in (0,\frac{bc}{d})$, then $f'(x)<0$; if
$x\in \left(\frac{bc}{d},+\infty\right)$, then $f'(x)>0$. It follows that
\begin{equation*}
\min\limits_{x\geq 0}\{f(x)\}=f\left(\frac{bc}{d}\right)=b^{p}d^{1-p}.
\end{equation*}

This completes the proof.\hfill$\Box$\vspace{.6pc}
\end{proof}

\subsection{\it The Minkowski inequality for dual quermassintegral sum of mixed intersection bodies}

In \cite{10}, Leng introduce the concept of $i$-{\it quermassintegral
difference function} of convex bodies as follows: If $K,D\in {\cal
K}^{n}$ and $D\subset K$, then $i$-quermassintegral difference
function of convex bodies $K$ and $D$, $D_{w_{i}}(K,D)$, is
defined by
\begin{equation*}
D_{w_{i}}(K,D)=W_{i}(K)-W_{i}(D)\qquad (0\leq i\leq n-1).
\end{equation*}

In the section, we first introduce a new concept, {\it dual
quermassintegral sum function}, as follows:

If $K,D\in \varphi^{n}$, then the dual quermassintegral sum
function of star bodies $K$ and $D$, $S_{\tilde{w}_{i}}(K,D)$, is
defined by
\begin{equation*}
S_{\tilde{w}_{i}}(K,D)=\tilde{W}_{i}(K)+\tilde{W}_{i}(D)\qquad (0\leq i\leq n-1).
\end{equation*}
When $i=0$, we have $S_{v}(K,D)=V(K)+V(D)$, which is called the
{\it dual volume sum function} of star bodies $K$ and $L$.

The following Minkowski inequality for mixed intersection bodies
will be established: If $K, L\in \varphi^{n}$, and $0\leq i< n$
and $0<j<n-1$, then
\setcounter{equation}{-1}
\begin{equation}
\tilde{W}_{i}(I_{j}(K,L))^{n-1}\leq \tilde{W}_{i}(I\!K)^{n-j-1}
\tilde{W}_{i}(I\!L)^{j},
\end{equation}
with equality if and only if $K$ and $L$ are dilates.\pagebreak

This is just the special case $D=D'=\O$ of the following inequality.

\begin{theor}[\!]
If $K, L, D, D'\in \varphi^{n}$. Let $D'$ is a dilate copy of $D${\rm ,} and $0\leq i<
n$ and $0<j<n-1${\rm ,} then
\begin{equation}
S_{\tilde{w}_{i}}(I_{j}(K,L),I_{j}(D,D'))^{n-1}\leq S_{\tilde{w}_{i}}
(I\!K, I\!D)^{n-j-1}S_{\tilde{w}_{i}}(I\!L, I\!D')^{j},
\end{equation}
with equality if and only if $K$ and $L$ are dilates.
\end{theor}

\begin{proof}
In view of the special case of Lemma B, we obtain that
\begin{equation}
\tilde{v}_{j}(K\cap E_{u},L\cap E_{u})^{n-i}\leq v(K\cap E_{u})^{\frac{(n-i)(n-j-1)}{n-1}}v(L\cap E_{u})^{\frac{j(n-i)}{n-1}}
\end{equation}
with equality if and only if $K\cap E_{u}$ and $L\cap E_{u}$ are
dilates. It follows if and only if $K$ and $L$ are
dilates \cite{20}.

From Lemma A, eq. (2.1.2) and in view of Minkowski inequality for
integral \cite{9}, we have for $i<n-1$,
\begin{align}
n\tilde{W}_{i}(I_{j}(K,L)) &=\left(\|\tilde{v}_{j}(K\cap E_{u},L\cap E_{u})\|_{n-i}\right)^{n-i}\nonumber\\[.5pc]
&\leq \left(\|v(K\cap E_{u})^{\frac{n-j-1}{n-1}}v(L\cap E_{u})^{\frac{j}{n-1}}\|_{n-i}\right)^{n-i}\nonumber\\[.5pc]
&\leq \left(\|v(K\cap E_{u}\|_{n-i}\right)^{\frac{(n-i)(n-j-1)}{n-1}}\left(\|v(K\cap E_{u}\|_{n-i}\right)^{\frac{j(n-i)}{n-1}}\nonumber\\[.5pc]
&=(n\tilde{W}_{i}(I\!K))^{\frac{(n-j-1)}{n-1}}(n\tilde{W}_{i}(I\!L))^{\frac{j}{n-1}}\nonumber\\[.5pc]
&=n\tilde{W}_{i}(I\!K)^{\frac{(n-j-1)}{n-1}}\tilde{W}_{i}(I\!L)^{\frac{j}{n-1}}.
\end{align}

In view of the conditions of (2.1.2) and Minkowski inequality for
integral, it follows that the equality holds if and only if $K$
and $L$ are dilates.

Moreover, we consider the case of $i=n-1$ of the inequality
(2.1.3). If $i=n-1$, inequality (2.1.3) reduces to
\begin{equation*}
\tilde{W}_{n-1}(I_{j}(K,L))^{n-1}\leq \tilde{W}_{n-1}(I\!K)^{n-j-1}\tilde{W}_{n-1}(I\!L)^{j}.\tag{*}
\end{equation*}
From Lemma A, (*) changes to
\begin{align*}
&\left(\int_{S^{n-1}}\tilde{v}_{j}(K\cap E_{u},L\cap E_{u})\hbox{d}S(u)\right)^{n-1}\\[.5pc]
&\quad\, \leq \left(\int_{S^{n-1}}v(K\cap E_{u})
\hbox{d}S(u)\right)^{n-j-1}\left(\int_{S^{n-1}}v(L\cap
E_{u})\hbox{d}S(u)\right)^{j}.\tag{**}
\end{align*}
On the other hand, integrating
both sides of (2.1.2) and in view of H\"{o}lder inequality for
integral, we obtain
\begin{align*}
&\int_{S^{n-1}}\tilde{v}_{j}(K\cap E_{u},L\cap E_{u})\hbox{d}S(u)\\[.5pc]
&\qquad\, \leq
\int_{S^{n-1}}v(K\cap E_{u})^{\frac{n-j-1}{n-1}}v(L\cap E_{u})^{\frac{j}{n-1}}\hbox{d}S(u)\\[.5pc]
&\qquad\, \leq  \left(\int_{S^{n-1}}v(K\cap
E_{u})\hbox{d}S(u)\right)^{\frac{n-j-1}{n-1}}\left(\int_{S^{n-1}}v(L\cap
E_{u})\hbox{d}S(u)\right)^{\frac{j}{n-1}}.
\end{align*}

Moreover, from inequality (2.1.3), we obtain
\begin{equation*}
\tilde{W}_{i}(I_{j}(K,L))^{n-1}\leq \tilde{W}_{i}(I\!K)^{n-j-1}\tilde{W}_{i}(I\!L)^{j},
\end{equation*}
with equality if and only if $K$ and $L$ are dilates, and
\begin{equation*}
\tilde{W}_{i}(I_{j}(D,D'))^{n-1}=\tilde{W}_{i}(I\!D)^{n-j-1}
\tilde{W}_{i}(I\!D')^{j}.
\end{equation*}
Hence, from the inequality in Lemma C, we have
\begin{align*}
S_{\tilde{w}_{i}}(I_{j}(K,L),I_{j}(D,D')) &\leq \tilde{W}_{i}(I\!K)^{(n-j-1)/(n-1)}
\tilde{W}_{i}(I\!L)^{j/(n-1)}\\[.3pc]
 &\quad\, +\tilde{W}_{i}(I\!D)^{(n-j-1)/(n-1)}\tilde{W}_{i}(I\!D')^{j/(n-1)}\\[.3pc]
&\leq S_{\tilde{w}}(I\!K, I\!D)^{n-j-1}S_{\tilde{w}}(I\!L, I\!D')^{j}.
\end{align*}
The proof of Theorem 2.1.1 is complete.\hfill$\Box$
\end{proof}

\setcounter{subsubsection}{0}
\begin{rem}
{\rm Taking $D=D'=\O$ and $j=1$ to (2.1.1), (2.1.1) changes to
\begin{equation*}
\tilde{W}_{i}(I_{1}(K,L))^{n-1}\leq \tilde{W}_{i}(I\!K)^{n-2}\tilde{W}_{i}
(I\!L),
\end{equation*}
with equality if and only if $K$ and $L$ are dilates.

This is just a dual form of the following inequality which was
given by Lutwak \cite{15}.}
\end{rem}

{\it The Minkowski inequality for mixed projection bodies}. If
$K, L\in {\cal K}^{n}$, and $0\leq i< n$, then
\begin{equation*}
W_{i}(\Pi_{1}(K,L))^{n-1}\geq W_{i}(\Pi K)^{n-2}W_{i}(\Pi L),
\end{equation*}
with equality if and only if $K$ and $L$ are homothetic.

A somewhat surprising consequence of Theorem 2.1.1 is the
following version.

\begin{theor}[\!]
If $K,L\in \eta\subset\varphi^{n}${\rm ,}
and $0\leq i<n${\rm ,} while $0<j<n-1$ and if either
\begin{equation}
\tilde{W}_{i}(I_{j}(K,M))=\tilde{W}_{i}(I_{j}(L,M)),\quad \hbox{for all}\ M\in \eta
\end{equation}
or
\begin{equation}
\tilde{W}_{i}(I_{j}(M,K))=\tilde{W}_{i}(I_{j}(M,L)),\quad \hbox{for all}\ M\in \eta
\end{equation}
hold{\rm ,} then it follows that $K=L${\rm ,} up to translation.
\end{theor}

\begin{proof}
Suppose (2.1.4) holds. Take $K$ for $M$ in
(2.1.4) and use the inequality (2.1.1). We obtain
\begin{equation*}
\tilde{W}_{i}(I\!K)=\tilde{W}_{i}(I_{j}(L,K))\leq \tilde{W}_{i}(I\!L)^{\frac{(n-j-1)}{n-1}}
\tilde{W}_{i}(I\!K)^{\frac{j}{n-1}}
\end{equation*}
with equality if and only if $K$ is a dilation of $L$.

Hence
\begin{equation*}
\tilde{W}_{i}(I\!K)\leq\tilde{W}_{i}(I\!L)
\end{equation*}
with equality if and only if $K$ is a dilation of $L$.

Similarly, take $L$ for $M$ in (2.1.4) and use again the
inequality (2.1.1). We get
\begin{equation*}
\tilde{W}_{i}(I\!K)\geq \tilde{W}_{i}(I\!L),
\end{equation*}
with equality if and only if $K$ is a dilation of $L$.

Hence
\begin{equation*}
\tilde{W}_{i}(I\!K)= \tilde{W}_{i}(I\!L)
\end{equation*}
and $K$ is a dilation of $L$. In view of the fact that the intersection bodies are
centered, there~exist $\lambda>0$ such that $K=\lambda L$,
and $\lambda^{(n-1)(n-i)}=1$, for $0\leq i<n-1$. Therefore
$\lambda=1$.~\hfill$\Box$
\end{proof}

Similar sort of argument shows that condition (2.1.5)
implies that $K$ and $L$ must be translates.

\setcounter{subsubsection}{1}
\begin{rem}
{\rm Theorem 2.1.2 is just the dual form of the
following `Theorem 5.4' which was given by Lutwak \cite{15}.}
\end{rem}

\setcounter{section}{5}
\setcounter{defin}{3}
\begin{theorem}[\!]
If $K,L\in \gamma\subset {\cal K}^{n}${\rm ,}
and $0\leq i<n${\rm ,} while $0<j<n-1$ and if either
\begin{equation*}
W_{i}(\Pi_{j}(K,M))=W_{i}(\Pi_{j}(L,M)),\quad \hbox{for}\ M\in \gamma
\end{equation*}
or
\begin{equation*}
W_{i}(\Pi_{j}(M,K))=W_{i}(\Pi_{j}(M,L)),\quad \hbox{for}\ M\in \gamma
\end{equation*}
hold{\rm ,} then it follows that $K=L${\rm ,} up to translation.
\end{theorem}

\setcounter{section}{2}
\subsection{\it The Aleksandrov--Fenchel inequality for mixed intersection bodies}

The Aleksandrov--Fenchel inequality for mixed intersection bodies
is as follows: If $K_{1},\ldots,K_{n-1}\in \varphi^{n}$, then
\begin{equation*}
V(I(K_{1},\ldots,K_{n-1}))\leq \prod_{j=1}^{r}V(I(\underbrace{K_{j},\ldots,K_{j}}_{r},K_{r+1},\ldots,K_{n-1}))
\end{equation*}
with equality if and only if $K_{1},\ldots,K_{n-1}$ are all
dilations of each other.

This is just the special case $i=0$ of the following.

\begin{theor}[\!]
If $K_{1},\ldots,K_{n-1}\in
\varphi^{n}${\rm ,} $0\leq i< n, 0<j<n-1$ and $0< r \leq n-1${\rm ,} then
\setcounter{equation}{0}
\begin{equation}
\tilde{W}_{i}(I(K_{1},\ldots,K_{n-1}))^{r}\leq\prod_{j=1}^{r}\tilde{W}_{i}(I(\underbrace{K_{j},\ldots,K_{j}}_{r},K_{r+1},\ldots,K_{n-1}))
\end{equation}
with equality if and only if
$K_{1},\ldots,K_{n-1}$ are all dilations of each other.
\end{theor}

\begin{proof}
When $i=1$, inequality (2.2.1) reduces to the inequality
in Lemma B. In the following, we suppose that $i<n-1$.\pagebreak

From (1.1.3) and (1.2.2), we have that
\begin{equation}
\tilde{W}_{i}(I(K_{1},\ldots,K_{n-1}))=\frac{1}{n}\int_{S^{n-1}}
\tilde{v}(K_{1}\cap E_{u},\ldots,K_{n-1}\cap E_{u})^{n-i}\hbox{d}S(u).
\end{equation}
By using the inequality in Lemma B, we easily get that
\begin{align}
&\tilde{v}(K_{1}\cap E_{u},\ldots,K_{n-1}\cap E_{u})^{n-i}\nonumber\\
&\quad\, \leq\left(\prod_{j=1}^{r}\tilde{v}(\underbrace{K_{j}\cap E_{u},\ldots,K_{j}\cap
E_{u}}_{r},K_{r+1}\cap E_{u},\ldots,K_{n-1}\cap
E_{u})\right)^{\frac{n-i}{r}},
\end{align}
with equality if and only if $K_{1}\cap E_{u},\ldots,K_{n-1}\cap E_{u}$ are all
dilations of each other. It follows if and only if
$K_{1},\ldots,K_{n-1}$ are all dilations of each other.

On the other hand, the H\"{o}lder's inequality can be stated as \cite{9}
\begin{equation}
\int_{S^{n-1}}\prod_{i=1}^{m}f_{i}(u)\hbox{d}S(u)\leq\prod_{i=1}^{m}\left(\int_{S^{n-1}}(f_{i}(u))^{m}\hbox{d}S(u)\right)^{1/m},
\end{equation}
with equality if and only if all $f_{i}$ are proportional.

From (2.2.2), (2.2.3) and (2.2.4), we obtain that

\begin{align*}
\tilde{W}_{i}(I(K_{1},\ldots,K_{n-1})) &=\frac{1}{n}\int_{S^{(n-1)}}\tilde{v}(K_{1}\cap E_{u},\ldots,K_{n-1}\cap E_{u})^{n-i}\hbox{d}S(u)\\
&\leq\frac{1}{n}\int_{S^{u}}\left(\prod_{j=1}^{r}\tilde{v}(\underbrace{K_{j}\cap E_{u},\ldots,K_{j}\cap E_{u}
}_{r},\right.\\
&\left.\quad\, \begin{array}{l} \\[2pc] \end{array}K_{r+1}\cap E_{u},\ldots,K_{n-1}\cap
E_{u})\right)^{\frac{n-i}{r}}\hbox{d}S(u)\\
&\leq \left(\prod_{j=1}^{r}\frac{1}{n}\int_{S^{n-1}}\tilde{v}(\underbrace{K_{j}\cap E_{u},\ldots,K_{j}\cap E_{u}
}_{r},\right.\\
&\quad\,\left. \begin{array}{l} \\[2pc] \end{array}K_{r+1}\cap E_{u},\ldots,K_{n-1}\cap
E_{u})^{n-i}\hbox{d}S(u)\right)^{r}\\
&=\left(\prod_{j=1}^{r}\tilde{W}_{i}(I\underbrace{(K_{j},\ldots,K_{j}}_{r},K_{r+1},\ldots,K_{n-1}))\right)^{r}.
\end{align*}

In view of the equality conditions (2.2.3) and (2.2.4), it
follows that the equality holds if and only if
$K_{1},\ldots,K_{n-1}$ are all dilations of each other.

The proof is complete.\hfill$\Box$
\end{proof}

\setcounter{subsubsection}{0}
\begin{rem}
{\rm The inequality (2.2.1) is just a dual form of
the following inequality which was given by Lutwak \cite{15}.}
\end{rem}

\noindent {\it The Aleksandrov--Fenchel inequality for mixed projection
bodies}. If $K_{1},\ldots,K_{n-1}\in {\cal K}^{n}$, $0\leq i<n,
1<j<n-1$ and $0< r \leq n-1$ then
\begin{equation*}
W_{i}(\Pi(K_{1},\ldots,K_{n-1}))^{r}\geq\prod_{j=1}^{r}W_{i}(\Pi(\underbrace{K_{j},\ldots,K_{j}}_{r},K_{r+1},\ldots,K_{n-1})).
\end{equation*}

From the case $r=n-1$ of inequality (2.2.1), it is as follows.

\setcounter{subsubsection}{0}
\begin{coro}$\left.\right.$\vspace{.5pc}

\noindent If $K_{1},\ldots,K_{n-1}\in
\varphi^{n}$, and $0\leq i< n$ then
\begin{equation}
\tilde{W}_{i}(I(K_{1},\ldots,K_{n-1}))^{n-1}\leq \tilde{W}_{i}
(I\!K_{1})\cdots\tilde{W}_{i}(I\!K_{n-1}),
\end{equation}
with equality if and only if $K_{1},\ldots,K_{n-1}$ are all
dilations of each other.
\end{coro}

\begin{rem}
{\rm Corollary (2.2.1) is similar to the following
`Theorem 5.2' which was given by Lutwak \cite{15}.}
\end{rem}

\setcounter{section}{5}
\setcounter{defin}{1}
\begin{theorem}[\!]
If $K_{1},\ldots,K_{n-1}\in {\cal K}^{n}$, and $0\leq i< n$ then
\begin{equation*}
W_{i}(\Pi(K_{1},\ldots,K_{n-1}))^{n-1}\geq W_{i}(\Pi
K_{1})\cdots W_{i}(\Pi K_{n-1}),
\end{equation*}
with equality if and only if $K_{1},\ldots,K_{n-1}$ are homothetic.
\end{theorem}

Taking $K_{1}=\cdots=K_{n-j-1}=K$ and $K_{n-j}=\cdots=K_{n-1}=L$
to (2.2.5), (2.2.5) reduces to (2.1.0). Taking
$K_{1}=\cdots=K_{r}=K$, $K_{r}=L$, and $K_{r+1}=\cdots=K_{n-1}=B$
to (2.2.1), (2.2.1) changes to the following.

\setcounter{section}{2}
\setcounter{subsubsection}{1}
\begin{coro}$\left.\right.$\vspace{.5pc}

\noindent If $K, L\in \varphi^{n}${\rm ,} and $0\leq i< n$ and $0\leq j<n-1${\rm ,} then
\begin{equation*}
\tilde{W}_{i}(I(\underbrace{K,\ldots,K}_{n-j-2},\underbrace{B,\ldots,B}_{j},L))^{n-j-1}\leq \tilde{W}_{i}(I_{j}K)^{n-j-2}\tilde{W}_{i}(I_{j}L),
\end{equation*}
with equality if and only if $K$ and $L$ are dilates.
\end{coro}

A somewhat surprising consequence of Corollary 2.2.2 is the
following version for mixed intersection bodies.

\setcounter{subsubsection}{1}
\begin{theor}[\!]
If $K,L\in \eta\subset \varphi^{n}${\rm ,}
$0\leq i<n-1, 0\leq j<n-1$ and if either
\begin{align}
&\tilde{W}_{i}(I(\underbrace{K,\ldots,K}_{n-j-1},\underbrace{B,\ldots,B}_{j},M))\nonumber\\[.5pc]
&\qquad\, =\tilde{W}_{i}(I(\underbrace{L,\ldots,L}_{n-j-1},\underbrace{B,\ldots,B}_{j},M)),\qquad
\hbox{for all}\ M\in \eta,
\end{align}
or
\begin{align}
&\tilde{W}_{i}(I(\underbrace{M,\ldots,M}_{n-j-1},\underbrace{B,\ldots,B}_{j},K))\nonumber\\[.5pc]
&\qquad\,=\tilde{W}_{i}(I(\underbrace{M,\ldots,M}_{n-j-1},\underbrace{B,\ldots,B}_{j},L)),\qquad\hbox{for all}\ M\in \eta,
\end{align}
hold{\rm ,} then it follows that $K$=$L${\rm ,} up to translation.
\end{theor}

\begin{proof}
Suppose (2.2.6) holds, take $K$ for $M$, use
Corollary~2.2.2, and get
\begin{align*}
\tilde{W}_{i}(I_{j}K)&=\tilde{W}_{i}(I(\underbrace{L,\ldots,L}_{n-j-1},\underbrace{B,\ldots,B}_{j},K))\\[.5pc]
&\leq \tilde{W}_{i}(I_{j}L)^{\frac{n-j-2}{n-j-1}}\tilde{W}_{i}(I_{j}K)^{\frac{1}{n-j-1}},
\end{align*}
with equality if and only if $K$ and $L$ are dilates.

Hence
\begin{equation*}
\tilde{W}_{i}(I_{j}K)\leq \tilde{W}_{i}(I_{j}L),
\end{equation*}
with equality if and only if $K$ and $L$ are dilates.

On the other hand, take $L$ for $M$, use Corollary 2.2.2 again,
and get
\begin{equation*}
\tilde{W}_{i}(I_{j}K)\geq \tilde{W}_{i}(I_{j}L),
\end{equation*}
with equality if and only if $K$ and $L$ are dilates.

Therefore
\begin{equation*}
\tilde{W}_{i}(I_{j}K)= \tilde{W}_{i}(I_{j}L),
\end{equation*}
where $K$ and $L$ are dilates and in view of the fact that the intersection bodies
are centered, there exists $\lambda>0$ such that $K=\lambda
L$. From (1.2.4), we have $\lambda^{(n-j-1)(n-i)}
\tilde{W}_{i}(I_{j}L)= \tilde{W}_{i}(I_{j}L),$ hence $\lambda=1$.

Similar argument shows that condition (2.2.7)
implies $K$=$L$, up to tanslation.\hfill$\Box$
\end{proof}

\begin{rem}
{\rm Taking $j=0$ to Theorem 2.2.2, it reduces to the
following:}\vspace{.6pc}

If $K,L\in \eta\subset \varphi^{n}$, $0\leq i<n$ and if either
\begin{equation*}
\tilde{W}_{i}(I_{1}(K,M))=\tilde{W}_{i}(I_{1}(L,M)),\qquad\hbox{for all}\ M\in \eta
\end{equation*}
or
\begin{equation*}
\tilde{W}_{i}(I_{1}(M,K))=\tilde{W}_{i}(I_{1}(M,L)),\qquad\hbox{for all}\ M\in \eta
\end{equation*}
hold, then it follows that $K=L$, up to translation.
\end{rem}

This is just the special case $j=1$ of Theorem~2.1.2.

\subsection{\it The Brunn--Minkowski inequality for mixed intersection bodies}

The Brunn--Minkowski inequality for intersection bodies, which will
be established is: If $K,L\in \varphi^{n}$, then
\begin{equation*}
V(I(K\tilde{+}L))^{1/n(n-1)}
\leq V(I\!K)^{1/n(n-1)}+V(I\!L)^{1/n(n-1)},
\end{equation*}
with equality if and only if $K$ and $L$ are dilates.

This is just the special case $i=0$ and $\alpha=1$ of the following.

\begin{theor}[\!]
If $K,L\in \varphi^{n}${\rm ,} and
$0\leq i<n${\rm ,} then for $0\leq \alpha \leq 1${\rm ,}
\setcounter{equation}{0}
\begin{align}
\tilde{W}_{i}(I(K\tilde{+}L))^{1/(n-i)(n-1)} &\leq
\tilde{W}_{i}(I(\alpha
K\tilde{+}(1-\alpha)L))^{1/(n-i)(n-1)}\nonumber\\
&\quad\, +\tilde{W}_{i}(I((1-\alpha)
K\tilde{+}\alpha L))^{1/(n-i)(n-1)},
\end{align}
with equality if and only if $(\alpha K\tilde{+}(1-\alpha)L)$ and
$(1-\alpha) K\tilde{+}\alpha L$ are dilates.
\end{theor}

\begin{proof}
Let $M=(L_{1},\ldots,L_{n-2})$, from (1.1.1), (1.1.3), (1.2.3)
and in view of the Minkowski inequality for
integral \cite{9}, we obtain that
\begin{align}
\tilde{W}_{i}(I(K\tilde{+}L,M))^{1/(n-i)} &=n^{-1/(n-i)}\|\rho(I(K\tilde{+}L,M),u)\|_{n-i}\nonumber\\[.3pc]
&=n^{-1/(n-i)}\|\rho(I(K,M)\tilde{+}I(L,M),u)\|_{n-i}\nonumber\\[.3pc]
&=n^{-1/(n-i)}\|\rho(I(K,M),u)+\rho(I(L,M),u)\|_{n-i}\nonumber\\[.3pc]
&\leq n^{-1/(n-i)}\left(\|\alpha\rho(I(K,M),u)\right.\nonumber\\
&\quad\,+(1-\alpha)\rho(I(L,M),u)\|_{n-i}+\|(1-\alpha)\nonumber\\
&\quad\,\left.\times\rho(I(K,M),u)+\alpha\rho(I(L,M),u)\|_{n-i}\right)\nonumber\\[.3pc]
&=n^{-1/(n-i)}\left(\|\rho(\alpha\cdot I(K,M)\tilde{+}(1-\alpha)\right.\nonumber\\
&\quad\, \times I(L,M),u)\|_{n-i}+\|\rho((1-\alpha)\cdot I(K,M)\nonumber\\
&\quad\,\left.\tilde{+}\alpha I(L,M),u)\|_{n-i}\right)\nonumber\\[.3pc]
&=n^{-1/(n-i)}\left(\|\rho(I(\alpha\cdot K\tilde{+}(1-\alpha)L,M),u)\|_{n-i}\right.\nonumber\\
&\quad\, \left.+\|\rho(I((1-\alpha)\cdot K\tilde{+}\alpha L,M),u)\|_{n-i}\right)\nonumber\\[.3pc]
&=\tilde{W}_{i}(I(\alpha\cdot K\tilde{+}(1-\alpha)L,M))^{1/(n-i)}\nonumber\\
&\quad\, +\tilde{W}_{i}(I(1-\alpha)\cdot K\tilde{+}\alpha L,M))^{1/(n-i)}.
\end{align}

On the other hand, taking $L_{1}=\cdots=L_{n-2}=K\tilde{+}L$ to
(2.3.2) and apply the inequality (2.1.0) twice, we get
\begin{align}
\tilde{W}_{i}(I(K\tilde{+}L))^{1/(n-i)}&\leq
\tilde{W}_{i}(I_{n-2}(\alpha\cdot
K\tilde{+}(1-\alpha)L,K\tilde{+}L))^{1/(n-i)}\nonumber\\
&\quad\, +\tilde{W}_{i}(I_{n-2}((1-\alpha)\cdot
K\tilde{+}\alpha L,K\tilde{+}L))^{1/(n-i)}\nonumber\\[.5pc]
&\leq \tilde{W}_{i}(I(\alpha\cdot
K\tilde{+}(1-\alpha)L))^{1/(n-1)(n-i)}\nonumber\\
&\quad\, \times\tilde{W}_{i}(I(K\tilde{+}L))^{(n-2)/(n-1)(n-i)}\nonumber\\
&\quad\,+\tilde{W}_{i}(I((1-\alpha)\cdot K\tilde{+}\alpha
L))^{1/(n-1)(n-i)}\nonumber\\
&\quad\, \times \tilde{W}_{i}(I(K\tilde{+}L))^{(n-2)/(n-1)(n-i)},
\end{align}
with equality if and only if $\alpha\cdot
K\tilde{+}(1-\alpha)L$, $(1-\alpha)\cdot K\tilde{+}\alpha L$ and
$M=K\tilde{+}L$ are dilates. Combining this with the equality
condition of (2.3.2), it follows that the condition holds if and
only if $K$ and $L$ are dilates.

Dividing both sides of (2.3.3) by
$\tilde{W}_{i}(I(K\tilde{+}L))^{(n-2)/(n-1)(n-i)}$, we get the
inequality (2.3.1).

The proof is complete.\hfill$\Box$
\end{proof}

Taking $\alpha=1$ in inequality (2.3.1), we have the following.

\setcounter{subsubsection}{0}
\begin{coro}$\left.\right.$\vspace{.5pc}

\noindent If $K,L\in \varphi^{n}${\rm ,} and $0\leq i<n${\rm ,} then
\begin{equation}
\tilde{W}_{i}(I(K\tilde{+}L))^{1/(n-i)(n-1)}\leq
\tilde{W}_{i}(I\!K)^{1/(n-i)(n-1)}+\tilde{W}_{i}(I\!L)^{1/(n-i)(n-1)},
\end{equation}
with equality if and only if $K$ and $L$ are dilates.
\end{coro}

\setcounter{subsubsection}{0}
\begin{rem}
{\rm From inequalities (2.3.4) and (1.2.4), we obtain that
\begin{align}
\tilde{W}_{i}(I(\alpha
K\tilde{+}(1-\alpha)L))^{1/(n-i)(n-1)}&\leq \alpha \tilde{W}_{i}(I\!K)^{1/(n-i)(n-1)}\nonumber\\
&\quad\, +(1-\alpha)\tilde{W}_{i}(I\!L)^{1/(n-i)(n-1)}
\end{align}
and
\begin{align}
\tilde{W}_{i}(I((1-\alpha)
K\tilde{+}\alpha L))^{1/(n-i)(n-1)}&\leq (1-\alpha) \tilde{W}_{i}(I\!K )^{1/(n-i)(n-1)}\nonumber\\
&\quad\, +\alpha\tilde{W}_{i}(I\!L)^{1/(n-i)(n-1)}.
\end{align}

From (2.3.5) and (2.3.6), we obtain that
\begin{align*}
&\tilde{W}_{i}(I(\alpha
K\tilde{+}(1-\alpha)L))^{1/(n-i)(n-1)}+\tilde{W}_{i}(I((1-\alpha)
K\tilde{+}\alpha L))^{1/(n-i)(n-1)}\\
&\quad\, \leq \tilde{W}_{i}(I\!K)^{1/(n-i)(n-1)}
+\tilde{W}_{i}(I\!L)^{1/(n-i)(n-1)}.
\end{align*}

This shows that inequality (2.3.1) is a strengthened form of
inequality (2.3.4).}
\end{rem}

\begin{rem}
{\rm Inequality (2.3.4) is just a dual form of
the following inequality which was given by Lutwak \cite{15}.}
\end{rem}

\noindent{\it The Brunn--Minkowski inequality for mixed projection bodies.}
If $K,L\in {\cal K}^{n}$, and $0\leq i<n,$ then
\begin{equation*}
W_{i}(\Pi(K+L))^{1/(n-i)(n-1)}\geq
W_{i}(\Pi K)^{1/(n-i)(n-1)}+W_{i}(\Pi L)^{1/(n-i)(n-1)},
\end{equation*}
with equality if and only if $K$ and $L$ are homothetic.

\section*{Acknowledgements}

This work was supported by the National Natural Sciences Foundation of China
(10271071), and by Academic Mainstay of Middle-age and Youth Foundation
of Shandong Province of China.

\end{document}